\documentclass[11pt]{article}
\usepackage{amssymb}

\usepackage{latexsym}
\usepackage{amsfonts}
\usepackage{amsthm}
\usepackage{amsmath}
\usepackage{latexcad}
\newtheorem{defin}{Definition}[section]
\newtheorem{theo}[defin]{Theorem}
\newtheorem{rem}[defin]{Remark}

\newtheorem{prop}[defin]{Proposition}
\newtheorem{cor}[defin]{Corollary}

\begin{document}
\date{February 2002}
\title{\Large Combinatorics and Topology of partitions \\
of spherical measures by 2 and 3 fans}

\author{{\Large Rade T. \v Zivaljevi\' c}  \\
       {Mathematics Institute SANU, Belgrade}}
\maketitle

\begin{abstract}
An arrangement of $k$-semilines in the Euclidean (projective)
plane or on the $2$-sphere is called a $k$-fan if all semilines
start from the same point. A $k$-fan is an $\alpha$-partition for
a probability measure $\mu$ if $\mu(\sigma_i)=\alpha_i$ for each
$i=1,\ldots,k$ where $\{\sigma_i\}_{i=1}^k$ are conical sectors
associated with the $k$-fan and $\alpha = (\alpha_1,\ldots
,\alpha_k)$.  The problem whether for a given collection of
measures $\mu_1,\ldots ,\mu_m$ and given $\alpha=(\alpha_1,\ldots
,\alpha_k)$ there exists a simultaneous $\alpha$-partition by a
$k$-fan was raised and studied in \cite{BaMa2001} in connection
with some partition problems in Discrete and Computational
Geometry\footnote{See the references \cite{Aki2000},
\cite{Bes2000}, \cite{Ito2000}, \cite{Kan1999}, \cite{Sak1998} and
\cite{VreZiv2001}.}. The set of all $\alpha = (\alpha_1,\ldots
,\alpha_m)$ such that for any collection of probability measures
$\mu_1,\ldots ,\mu_m$ there exists a common $\alpha$-partition by
a $k$-fan is denoted by ${\cal A}_{m,k}$. It was shown in
\cite{BaMa2001} that the interesting cases of the problem are
$(m,k) = (3,2),(2,3),(2,4)$. We prove, as a central result of this
paper, that ${\cal A}_{3,2} = \{(s,t)\in \mathbb{R}^2\mid s+t=1
\mbox{ {\rm and} } s,t>0\}$. The result follows from the fact that
under mild conditions there does not exist a $Q_{4n}$-equivariant
map $f : S^3\rightarrow V\setminus {\cal A}(\alpha)$ where ${\cal
A}(\alpha)$ is a $Q_{4n}$-invariant, linear subspace arrangement
in a $Q_{4n}$-representation $V$, where $Q_{4n}$ is the
generalized quaternion group. This fact is established by showing
that an appropriate obstruction in the group $\Omega_1(Q_{4n})$ of
$Q_{4n}$-bordisms does not vanish.
\end{abstract}

\section{Introduction}

 A $k$-fan ${\mathfrak p} = (x;l_1,l_2,\ldots ,l_k)$ on the sphere
$S^2$ consists of a point $x$, called the center of the fan, and
$k$ great semicircles $l_1,\ldots ,l_k$ emanating from $x$. It is
always assumed that semicircles in a $k$-fan ${\mathfrak
p}=(x;l_1,\ldots ,l_k)$ are enumerated counter clockwise, in an
agreement with the standard orientation on the ambient $2$-sphere.
Sometimes it is more convenient to use the notation ${\mathfrak p}
= (x; \sigma_1,\sigma_2,\ldots,\sigma_k)$, where $\sigma_i$ is the
open angular sector between $l_i$ and $l_{i+1}, \, i=1,\ldots ,k$.
Here we adopt the circular order $1\prec 2\prec\ldots\prec k\prec
1$ of indices and their addition ``modulo $k$'', so for example
$\sigma_k$ is the open angular sector between $l_k$ and
$l_1=l_{k+1}$.

Let $\mu_1,\mu_2,\ldots,\mu_m$ be Borel probability measures on
$S^2$. We assume in this paper that all measures $\mu_j$ are {\em
proper} in the sense that $\mu_j([a,b]) =0$ for any circular arc
$[a,b]\subset S^2$ and that $\mu_j(U)>0$ for each nonempty open
set $U\subset S^2$. All the results, appropriately reformulated,
can be extended by a standard limit argument to more general
measures, including the counting measures of finite sets, see
\cite{BaMa2001}, \cite{elisat}, \cite{T-V} for related examples.

 Let $(\alpha_1,\alpha_2,\ldots
,\alpha_k)$ be a vector of positive real numbers where
$\alpha_1+\alpha_2 +\ldots +\alpha_k = 1$. Following
\cite{BaMa2001}, but taking into account our simplifying
assumptions that we deal only with {\em proper} measures, we say
that a $k$-fan $(x;l_1,\ldots ,l_k)$ together with the associated
conical sectors $\sigma_1,\ldots ,\sigma_k$ is an
$\alpha$-partition for the collection $\{\mu_j\}_{j=1}^m$ of
measures if $\mu_j(\sigma_i)=\alpha_i$ for all $i=1,\ldots ,k$ and
$j=1,\ldots,m$. Note that for more general measures, say for the
weak limits of proper measures, instead of the equality
$\nu(\sigma_i)=\alpha_i$ one could use a pair of inequalities
$\nu(\sigma_i)\leq\alpha_i\leq\nu(\bar{\sigma_i})$.

\begin{defin}
\label{def:admissible} A vector $\alpha\in \mathbb{R}^k$ is called
admissible, or more precisely $(m,k)$-{\em admissible}, if for any
collection of $m$ (proper) measures on $S^2$, there exists a
simultaneous $\alpha$-partition. The collection of all
$(m,k)$-admissible vectors is denoted by ${\cal A}_{m,k}$.
\end{defin}

\medskip\noindent
{\bf Problem} (\cite{BaMa2001}): Find a characterization of the
set ${\cal A}_{m,k}$ or at least describe some of its nontrivial
properties. Equivalently, the problem is to decide for what
combinations of $m,k$ and $\alpha\in \mathbb{R}^k$ one can
guarantee that for any collection of measures ${\cal M} =
\{\mu_1,\mu_2,\ldots ,\mu_m\}$, there exist an $\alpha$-partition
for ${\cal M}$.

\medskip

The analysis given in \cite{BaMa2001} shows that the most
interesting cases for the problem of the existence of
$\alpha$-partitions are $(m,k) = (3,2), (2,3), (2,4)$. In the case
$(m,k)=(3,2)$ it was shown that $\{(\frac{1}{2},\frac{1}{2}),
(\frac{1}{3},\frac{2}{3}), (\frac{2}{3},\frac{1}{3})\}\subseteq
{\cal A}_{3,2}.$

\medskip
In this paper we prove, Theorem~\ref{thm:prva}, that
 \[{\cal A}_{3,2}=\{(s,t)\in \mathbb{R}^2 \mid s+t =1 { \mbox{ \rm and } } s,t >
0\}.\]

\section{The candidate space/test map paradigm}
\label{sec:paradigm}

Imre B\' ar\' any and Ji\v ri Matou\v sek demonstrated in
\cite{BaMa2001} that the problem of $\alpha$-partitioning of
$m$-measures on $S^2$ by a spherical $k$-fan can be in a very
elegant way reduced to the problem of the existence of equivariant
maps. B\' ar\' any and Matou\v sek introduced, for a given
positive vector $\alpha = (\frac{a_1}{n}, \frac{a_2}{n}, \ldots ,
\frac{a_k}{n})\in\frac{1}{n}\,\mathbb{Z}^k\subset \mathbb{Q}^k$, a
$D_{2n}$-invariant, linear subspace arrangement ${\cal A} = {\cal
A}(\alpha)={\cal A}_m^k(\alpha)$ of $(W_n)^{\oplus(m-1)}$, where
$W_n :=\{x\in \mathbb{R}^n \mid x_1+\ldots +x_n =0\}$ and $D_{2n}$
is the dihedral group of order $2n$. They showed that an
$\alpha$-partition exists {\em if} there does not exist a
$D_{2n}$-equivariant map $f : V_2(\mathbb{R}^3)\rightarrow
(W_n)^{\oplus(m-1)}\setminus\cup{\cal A}$, from the Stiefel
manifold $V_2(\mathbb{R}^3)\cong SO(3)$ of all orthonormal
$2$-frames in $\mathbb{R}^3$ to the complement $M({\cal A}):=
(W_n)^{\oplus(m-1)}\setminus\cup{\cal A}$ of the arrangement
${\cal A}$.

\medskip
Here is a brief outline of this construction. One starts with the
definition of the {\em candidate} or {\em configuration space}
$X_{\mu_1}$ which is constructed with the aid of the first measure
$\mu_1$. By definition ${\mathfrak q} = (x;l_1,\ldots ,l_n)\in
X_{\mu_1}$ if and only if $\mu_1(\sigma_i)=\frac{1}{n}$ for each
$i=1,\ldots ,n$. In other words $X_{\mu_1}$ is the collection of
all $n$-fans on the sphere which form an equipartition of the
first measure $\mu_1$. Note that ${\mathfrak q} = (x;l_1,\ldots
,l_n)\in X_{\mu_1}$ is uniquely determined by the pair $(x,l_1)$
or equivalently the pair $(x,y)$, where $y$ is the unit tangent
vector to $l_1$ at $x$. Hence, the space $X_{\mu_1}$ is
homeomorphic to the Stiefel manifold $V_2(\mathbb{R}^3)\cong
SO(3)$ of all orthonormal $2$-frames, respectively all
orthonormal, positive $3$-frames in $\mathbb{R}^3$.

In order to check how far is a $n$-fan ${\mathfrak q}\in
X_{\mu_1}$ from being an equipartition for the measure $\mu_j, \,
j\geq 2$, one introduces a {\em test map} $f_j :
X_{\mu_1}\rightarrow \mathbb{R}^n$ by
 \[ f_j({\mathfrak q}) = (\mu_j(\sigma_1)-\textstyle\frac{1}{n},\ldots
,\mu_j(\sigma_n)-\textstyle\frac{1}{n}), \] where as before
$\sigma_i$ is the sector on the sphere bounded by $l_i$ and
$l_{i+1}$. Taken together, the maps $\{f_j\}_{j=2}^m$ define a
test map $F: X_{\mu_1}\rightarrow (\mathbb{R}^n)^{\oplus(m-1)}$
where
 \[ F(\mathfrak q) := (f_2({\mathfrak q}),\ldots ,f_m({\mathfrak q}))\in
 (\mathbb{R}^n)^{\oplus (m-1)}.\]
Sometimes it is convenient to interpret the target space
$(\mathbb{R}^n)^{\oplus (m-1)}$ as the space ${\rm
Mat}_{(m-1)\times n}(\mathbb{R})$ of all $(m-1)$ by $n$ matrices.
Note also that ${\rm Im}(f_j)\subset W_n := \{x\in \mathbb{R}^n
\mid x_1+\ldots +x_n=0\}$ which implies that the matrix
$F({\mathfrak q})$ has an additional property that the entries in
each row add up to zero. The column decomposition ${\rm
Mat}_{(m-1)\times n}({\mathbb R})\cong ({\mathbb
R}^{(m-1)})^{\oplus n} \cong L_1\oplus\ldots\oplus L_n, \, x =
x_1\oplus\ldots\oplus x_n = (x_1,\ldots, x_n)$, where $L_i\cong
{\mathbb R}^{m-1}$ is also useful. All these target spaces are
identified, so the actual {\em test space} $W = W_n^{\oplus
(m-1)}$ is seen as a vector subspace in each of spaces
\[
{\rm Mat} := {\rm Mat}_{(m-1)\times n}({\mathbb R})\cong
({\mathbb R}^{(m-1)})^{\oplus n}
\cong ({\mathbb R}^{n})^{\oplus (m-1)}.
\]

The dihedral group $D_{2n}$, defined as the group generated by $E$ and $J$
subject to the relations
\begin{equation}
\label{eqn:dihedral} E^n=1  \qquad JEJ = E^{n-1},
\end{equation}
acts on $\mathbb{R}^n, W_n$ and similarly on linear spaces ${\rm
Mat}_{(m-1)\times n}({\mathbb R})$ and $W_n^{\oplus(m-1)}$ by
\[ E(x_1,x_2,\ldots ,x_n) = (x_2,\ldots, x_n, x_1),\,
J(x_1,x_2,\ldots , x_n) = (x_n, x_{n-1},\ldots ,x_1).\] The group
$D_{2n}$ acts also on the configuration space $X_{\mu_1}$ by
\[ E(x;l_1,l_2,\ldots ,l_n) = (x;l_2,\ldots ,l_n,l_1), \,
J(x;l_1,l_2,\ldots ,l_n) = (-x;l_1,l_n,\ldots ,l_2).\] Perhaps the
action of $J$ appears more natural if expressed as
\[J(x;\sigma_1,\sigma_2,\ldots,\sigma_n) =
(-x;\sigma_n,\sigma_{n-1},\ldots, \sigma_1)\] where $\sigma_i$ are
the sectors on $S^2$ associated to the $n$-fan $(x;l_1,\ldots,
l_n)$. The following proposition doesn't require a proof.

\begin{prop}
The test map $F : X_{\mu_1}\rightarrow W_n^{\oplus(m-1)}$ is
$D_{2n}$-equivariant.
\end{prop}

\medskip
Let us associate to each vector $\alpha = (\frac{a_1}{n},\ldots
,\frac{a_k}{n})$ a space $L=L(\alpha)\subset
(W_n)^{\oplus(m-1)}\subset {\rm Mat}$,
\begin{equation}
\label{eqn:L}
 L = L(\alpha) := \{x\in (W_n)^{\oplus{(m-1)}}\mid
 z_1(x)=z_2(x)= \ldots =z_k(x) = 0\},
\end{equation}
where $z_i : {\rm Mat}\rightarrow W_n$ are the linear maps defined
by
\begin{equation}
\label{eqn:formsforms}
\begin{array}{ccl}
z_1 & = & x_1 +  \ldots + x_{a_1} \\
z_2 & = & x_{a_1+1} +\ldots + x_{a_1+a_2}\\
z_3 & = & x_{a_1+a_{2}+1}+\ldots + x_{a_1+a_2+a_3}\\
&&\ldots\ldots\ldots\ldots \\
 z_k & = & x_{a_1+\ldots +
a_{k-1}+1}+\ldots + x_{n}
\end{array}
\end{equation}
and $x_j : {\rm Mat}\rightarrow L_j$ are projections on the
corresponding column spaces.

Let ${\cal A} = {\cal A}(\alpha)$ be the smallest
$D_{2n}$-invariant, linear subspace arrangement in
$(W_n)^{\oplus(m-1)}$ which contains $L(\alpha)$. Hence, ${\cal
A}(\alpha)$ is the arrangement which contains all subspaces of the
form $L_g(\alpha):= g(L(\alpha)), \, g\in D_{2n}$, together with
their intersections.

\begin{prop}
\label{prop:reduction} {\rm ({\cite{BaMa2001}})}
 Let $\alpha = (\frac{a_1}{n},
 \frac{a_2}{n},\ldots,\frac{a_k}{n})\in \mathbb{R}^k$ be a
 vector such that all $a_i$ are positive integers and
 $a_1+\ldots +a_k = n$. Let us suppose that there does not exist a
 $D_{2n}$-equivariant map $F : V_2(\mathbb{R}^3)\rightarrow
 M(\alpha)$, where $M(\alpha) := (W_n)^{\oplus(m-1)}\setminus \cup{\cal
 A}(\alpha)$. Then, $\alpha\in {\cal A}_{m,k}$, or in other
 words for any collection $\{\mu_j\}_{j=1}^m$ of (proper) measures on
 $S^2$, there always exists a $k$-fan which
 $\alpha$-partitions each of the measures $\mu_i$.
\end{prop}

\begin{rem}
{\rm The ``configuration space/test map'' scheme, used by B\' ar\'
any and Matou\v sek in the $\alpha$-partition problem, has been
successfully applied on numerous problems in discrete geometry and
combinatorics. The review papers \cite{Bar93}, \cite{Bjo91},
\cite{Mat94}, \cite{elisat}, \cite{guide1},\cite{guide2} provide
more information and give references to the original papers where
these ideas were introduced and developed. This method continues
to serve as a central guiding principle for applying topological
methods in discrete problems.}
\end{rem}

\section{The change of the group and equivariant maps}

If $G$ acts on a space $X$, say on the left, and
$H\vartriangleleft G$ is a normal subgroup, then the quotient
group $Q:=G/H$ acts on the orbit space $X/H$ by $gH(Hx):=H(gx)$.
For example if the cyclic group $G= \mathbb{Z}/{2n} \cong
\{1,\omega,\ldots ,\omega^{2n-1}\}$ acts freely on $S^3$, then
there is an induced action of $\mathbb{Z}/n\cong G/H=
(\mathbb{Z}/{2n})/(\mathbb{Z}/2)$ on $S^3/( \mathbb{Z}/2)\cong
RP^3$, where $H\cong \mathbb{Z}/2 = \{1,\omega^{n}\}\subset G$.

A $Q$-space $Z$ is also seen as a $G$-space via the obvious
homomorphism $\rho : G\rightarrow Q$ where for $g\in G$ and $z\in
Z, \, g\cdot z:=\rho(g)z$.

The following proposition is a topological analogue of the well
known ``extension of scalars'' equivalence from homological
algebra, \cite{Brown} Section~III.3.

\begin{prop}
\label{prop:homological} Suppose that $X$ and $Z$ are $G$-spaces
and $H$ is a normal subgroup of $G$ which acts trivially on $Z$.
There exists a $G$-equivariant map $\alpha : X\rightarrow Z$ if
and only if there exists a $Q$-equivariant map $\beta :
X/H\rightarrow Z$, where $Q:= G/H$ and $X/G$ and $Z/H = Z$ are
interpreted as $Q$-spaces.
\end{prop}

\medskip\noindent
{\bf Proof:} The obvious quotient map $p : X\rightarrow X/H$ is
$G$-equivariant. As before the $Q$-space $X/H$ is automatically a
$G$-space via the homomorphism $\rho : G\rightarrow Q$. If $\beta
:X/H\rightarrow Z$ is $Q$-equivariant, it is also $G$-equivariant
and the composition $\alpha:= p\circ \beta : X\rightarrow Z$ is
$G$-equivariant. Conversely, if $\alpha : X\rightarrow Z$ is a
$G$-equivariant map then, in light of the fact that $H$-acts
trivially on $Z$, there is a factorization $\alpha = p\circ\beta$
for some $\beta : X/H\rightarrow Z$, i.e. $\beta(Hx) = \alpha(x)$
for each $x\in X$. It is easily checked that $\beta$ is
$Q$-equivariant,

$$
\beta(gH\cdot Hx):=\alpha(gx)=g\alpha(x)=g\beta(Hx) =
(gH)\beta(Hx).\eqno{\square}
$$

\medskip
In this paper we pay a special attention to the class of {\em
generalized quaternion groups}, \cite{CaEi} p. 253, which are also
known as {\em binary dihedral groups}. Let $S^3 = S(\mathbb{H})=
Sp(1)$ be the group of all unit quaternions.  Let $\epsilon =
\epsilon_{2n} = \cos\frac{\pi}{n} + i\sin\frac{\pi}{n}\in
S(\mathbb{H})$ be a root of unity generating a subgroup of
$S(\mathbb{H})$ of order $2n$. Then
\[
Q_{4n} = \{1,\epsilon,\ldots ,\epsilon^{2n-1}, j, \epsilon j,
\ldots , \epsilon^{2n-1}j\}
\]
is a subgroup of $S^3$ of order $4n$ called {\em the generalized
quaternion group}. It is easily checked that $Q_{4n}$ is
isomorphic to the group freely generated by $\epsilon$ and $j$
subject to the relations

\[  \epsilon^n = j^2  \qquad \epsilon j\epsilon = j. \]

Let $H = \{1,\epsilon^n\}=\{1,-1\}\subset Q_{4n}$ and let $D_{2n}$
be the associated quotient group $D_{2n} = Q_{4n}/H$. Then
$D_{2n}$ is easily found to be isomorphic to the dihedral group of
order $2n$. Moreover, since $Q_{4n}/H\subset S^3/H\cong SO(3)$, we
identify $D_{2n}$ with a subgroup $D'_{2n}$ of $SO(3)$. The
following well known description of the homomorphism $\theta :
S^3\rightarrow SO(3)$ will help us give a more transparent picture
of the group $D'_{2n}$. Given $q\in S^3$, the isometry
$\theta(q)\in SO(3)$ is described by $\theta(q)x = qx{\bar q}$
where $x\in \mathbb{R}^3\cong {\rm Im}(\mathbb{H})$. If $q =
\cos\alpha + u\sin\alpha$, where $0\leq\alpha\leq\pi$, and $u\in
{\rm Im}(\mathbb{H})$, then $\theta(q)$ has a geometric
description as the positive rotation $R_u(2\alpha)$ around the
oriented axes determined by $u$, through the angle $2\alpha$.

\medskip
As a consequence we obtain that $\theta(\epsilon) =
R_i(\frac{2\pi}{n}), \theta(j)= R_j(\pi)$ and $D'_{2n}$ is
identified as a dihedral group of order $2n$. Note that an
abstract dihedral group is generated by two generators $E$ and $J$
subject to relations described in (\ref{eqn:dihedral}). We often
identify all these dihedral groups so in particular
\begin{equation}
\label{eqn:write}
\begin{array}{cccc}
E=R_i(\frac{2\pi}{n}), & J = R_j(\pi), & \theta(\epsilon) = E, &
\theta(j) = J.
\end{array}
\end{equation}

\medskip
The following simple and useful result allows us to change or
modify the action of the group $G$ on the domain $S^3$.

\begin{prop}
\label{prop:change} Suppose that $\gamma_i : G\times
S^3\rightarrow S^3, \gamma_i(g,x) = g\cdot_i x, \, i=1,2$ are two
actions of a finite group $G$ on the $3$-sphere $S^3$. Suppose
that the action $\gamma_1$ is free. Then there exists a
$G$-equivariant map $f : S^3\rightarrow S^3$ between these two
actions, that is a map such that for each $g\in G$ and $x\in S^3$,
\[    f(g\cdot_1 x) = g\cdot_2 f(x) \]
\end{prop}

\medskip\noindent
{\bf Proof:} The proof of this fact is routine and relies on the
fact that $S^3$ is a $2$-connected, $\gamma_1$-free, $CW$-complex
so there are no obstructions to extend equivariantly a map defined
on $0$-skeleton of $S^3$.\hfill$\square$

\medskip

\begin{cor}
\label{cor:change}
Suppose that $\gamma_i : G\times S^3\rightarrow
S^3, \gamma_i(g,x) = g\cdot_i x, \, i=1,2$ are two free actions of
a finite group $G$ on the $3$-sphere $S^3$ and let $Z$ be an
arbitrary $G$-space. Then there exists a $\gamma_1$-equivariant
map $f : S^3\rightarrow Z$ if and only if such a map exists for
the $\gamma_2$ action.
\end{cor}

As an illustration of the use of Proposition~\ref{prop:change} and
its Corollary~\ref{cor:change} we prove the following proposition
which is needed in Section~\ref{sec:2-fans}. Let us start with a
definition.
\begin{defin}
Suppose that the linear subspace arrangements ${\cal A}$ and
${\cal B}$ in $\mathbb{R}^n$ are both $G$-invariant for some
$G\subset GL(n,\mathbb{R})$. ${\cal A}$ and ${\cal B}$ are {\em
isomorphic} if there exists a non singular linear map $C:
\mathbb{R}^n\rightarrow \mathbb{R}^n$ such that $K\in {\cal A}
\Leftrightarrow C(K)\in {\cal B}$. If $C : \mathbb{R}^n\rightarrow
\mathbb{R}^n$ is $G$-equivariant, i.e. if $Cg = gC$ for each $g\in
G$, we say that ${\cal A}$ and ${\cal B}$ are $G$-isomorphic.
Finally ${\cal A}$ and ${\cal B}$ are {\em weakly} $G$-isomorphic
or $wG$-isomorphic if there exists an automorphism $\theta :
G\rightarrow G$ of the group $G$, such that $C(g\cdot x) =
\theta(g)\cdot C(x)$.
\end{defin}

\begin{prop}
\label{prop:spec}
If two arrangements ${\cal A}$ and ${\cal B}$
are $wG$-isomorphic, then there exists a $G$-equivariant map $f :
S^3\rightarrow M({\cal A})$ if and only if there exists a
$G$-equivariant map $g : S^3\rightarrow M({\cal B})$.
\end{prop}

Here is an example of two $wG$-isomorphic arrangements.

\begin{prop}
\label{prop:two} Let $\alpha = (\frac{1}{n},\frac{n-1}{n})$ and
$\beta = (\frac{p}{n},\frac{n-p}{n})$, where $p$ is an integer and
$1\leq p\leq n-1$. Then the arrangements ${\cal A}(\alpha)$ and
${\cal A}(\beta)$, interpreted as arrangements both in ${\rm
Mat}_{2\times n}\cong (\mathbb{R}^2)^{\oplus n}$ and its subspace
$W_n^{\oplus 2}$, are $wD_{2n}$-isomorphic.
\end{prop}

\medskip\noindent
{\bf Proof:} Let $x_i : (\mathbb{R}^2)^{\oplus n} \rightarrow
\mathbb{R}^2, \, i=1,\ldots ,n$ be the projection on the $i$-th
factor. Let $y_j : (\mathbb{R}^2)^{\oplus n} \rightarrow
\mathbb{R}^2, \, j=1,\ldots ,n$ be defined by $y_j =
x_{j+1}+\ldots +x_{j+p}$. The arrangement ${\cal A}(\alpha)$ is
generated by subspaces ${\rm Ker}(x_i)$ and similarly, ${\cal
A}(\beta)$ is generated by ${\rm Ker}(y_j)$. Then
\[
x_i\circ E= x_{i+1}, \quad y_i\circ E = y_{i+1}, \quad x_i\circ J
= x_{n-i+1}, \quad y_i\circ J = y_{n-i-p}
\]
and if $C:(\mathbb{R}^2)^{\oplus n}\rightarrow
(\mathbb{R}^2)^{\oplus n}$ is the linear map such that $y_i =
x_i\circ C$ then
\[ EC = CE \quad { \rm{ and }} \quad CJ = JE^{-p-1}C.\]
Hence, ${\cal A}(\alpha)$ and ${\cal B}(\beta)$ are
$wD_{2n}$-isomorphic with the associated automorphism $\theta :
D_{2n}\rightarrow D_{2n}$ defined by $\theta(E) = E$ and
$\theta(J) = JE^{-p-1}$. Note that if $p=2s-1$ is an odd integer
then $JE^{-p-1} = E^sJE^{-s}$ and $\theta$ is an inner
automorphism of $D_{2n}$. \hfill$\square$

\begin{cor} There exists a $D_{2n}$-equivariant map $f :
V_2(\mathbb{R}^3)\rightarrow W_n^{\oplus 2}\setminus\cup {\cal
A}(\alpha)$ if and only if there exists an equivariant map $g :
V_2(\mathbb{R}^3)\rightarrow W_n^{\oplus 2}\setminus\cup {\cal
A}(\beta)$

\end{cor}

\section{Evaluation of the group $\Omega_1(Q_{4n})$}
\label{sec:evaluation}

In this section we collect some elementary facts about the group
$\Omega_1(Q_{4n})$ of $Q_{4n}$-bordism classes of oriented,
$1$-dimensional $Q_{4n}$-manifolds. In light of natural
isomorphisms,
\begin{equation}
\label{eqn:natural} \Omega_1(G)\cong \Omega_1(BG)\cong H_1(G,
\mathbb{Z}),
\end{equation}
see \cite{CoFl} or \cite{Go} Lemma~2, the computation of
$\Omega_1(G)$ is reduced to evaluation of group homology. On the
other hand the {\em Abelization} functor $G\mapsto {\rm
Ab}(G):=G/[G,G]$ allows us to compute $H_1(X, \mathbb{Z})$, at
least for connected $X$, as the group ${\rm Ab}(\pi_1(X))$. In
particular we have isomorphisms
 \[
 H_1(G, \mathbb{Z}) = H_1(BG, \mathbb{Z})\cong {\rm Ab}(\pi_1(BG))
 \cong {\rm Ab}(G) .
 \]

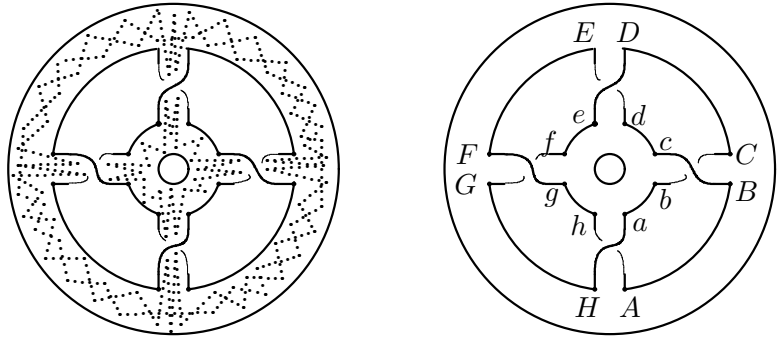
\placedrawing[htbp]{lim1.lp}{The class
$[\mathbb{Z}/4\times_{\mathbb{Z}/2}S^1]$ in $\Omega_1(
\mathbb{Z}/4)$.}{lim1}

\begin{prop}
\label{prop:mono} Suppose that $Q_{4n} = \{1,\epsilon,\ldots
,\epsilon^{2n-1}, j, \epsilon j, \ldots , \epsilon^{2n-1}j\}$ is a
generalized quaternion group of order $4n$ where $n$ is an odd
integer. Let $H = \{1,j,j^2,j^3\}$ be the subgroup generated by
$j$. Then the map,
$\Omega_1(H)\stackrel{I}{\rightarrow}\Omega_1(Q_{4n})$, induced by
the inclusion $H\stackrel{i}{\rightarrow} G$, is an isomorphism.
Hence, \[\Omega_1(Q_{4n})\cong \Omega_1(\mathbb{Z}/4) \cong
\mathbb{Z}/4.\]
\end{prop}

\medskip\noindent
{\bf Proof:} The commutativity of the following diagram, cf.
\cite{CoFl}, Section~III.20,
\begin{equation}
\label{eqn:square}
\begin{array}{ccc}
 \Omega_1(H) & \rightarrow & \Omega_1(Q_{4n}) \\
 \downarrow & & \downarrow \\
 H_1(H, \mathbb{Z}) & \rightarrow & H_1(Q_{4n}, \mathbb{Z})
\end{array}
\end{equation}
where the vertical arrows are isomorphisms, allows us to prove the
result on the level of homology. Let $A:=
\{\epsilon^{2j}\}_{j=1}^n$ be the subgroup of $Q_{4n}$ generated
by $\epsilon^2$. Since $A = [Q_{4n}, Q_{4n}]$ is the commutator
group, it is a normal subgroup of $Q_{4n}$. Moreover, since $n$ is
an odd number, $A\cap H=\{1\}$. As a consequence one obtains that
the following composition of homomorphisms
 \[
 H\rightarrow Q_{4n}\rightarrow
 Q_{4n}/A\stackrel{\cong}{\rightarrow} H
 \]
is an isomorphism. By applying the {\em Abelization} functor on
this sequence, we obtain that the horizontal, as well as the
vertical arrows in the diagram
\begin{equation}
\begin{array}{ccccc}
 H & \rightarrow & {\rm Ab}(Q_{4n}) &\rightarrow & H \\
 \downarrow & & \downarrow && \downarrow \\
 H_1(H, \mathbb{Z}) &\rightarrow &H_1(Q_{4n}, \mathbb{Z}) &
 \rightarrow & H_1(H, \mathbb{Z})
\end{array}
\end{equation}
are also isomorphisms. Hence, the map $H_1(H, \mathbb{Z})
\rightarrow H_1(Q_{4n})$ is an isomorphism.\hfill$\square$

\medskip
The following proposition is obtained by a similar proof.
\begin{prop}
\label{prop:basic} The map $\Omega_1(\mathbb{Z}/2)\rightarrow
\Omega_1(\mathbb{Z}/4)$, induced by the map $\omega\mapsto j^2$,
where $\mathbb{Z}/2 = \{1,\omega\}$ and $\mathbb{Z}/4 =
\{1,j,j^2,j^3\}$, is a monomorphism. Moreover, the image of the
generator of $\Omega(\mathbb{Z}/2)\cong H_1(\mathbb{Z}/2,
\mathbb{Z}) = \mathbb{Z}/2$ in $\Omega_1(\mathbb{Z}/4)\cong
\mathbb{Z}/4$ is the $1$-dimensional, $\mathbb{Z}/4$-manifold
$\mathbb{Z}/4\times_{\mathbb{Z}/2}S^1$, where $S^1$ is seen as a
$\mathbb{Z}/2$-manifold.
\end{prop}

\medskip
Figure~\ref{lim1} gives a pictorial proof of
Proposition~\ref{prop:basic}. It shows that the union of circles
$C_1 = ``{\rm ABcdEFgh}''$ and $C_2 = ``{\rm abCDefGH}''$ is
$\mathbb{Z}/4$-bordant to the union of two (positive)
$\mathbb{Z}/4$-circles, hence it represents the element $2\in
\mathbb{Z}/4\cong\Omega_1(\mathbb{Z}/4)$.

\begin{cor}
\label{cor:vazan} Let $\Delta := Q_{4n}\times_{\mathbb{Z}/2}S^1 =
Q_{4n}\times_{\mathbb{Z}/4}M^1$ where $\mathbb{Z}/4$ and
$\mathbb{Z}/2$ are respectively the subgroups of $Q_{4n}$
generated by $j$ and $j^2$, the action of $j^2$ on $S^1$ is
antipodal and $M^1 := \mathbb{Z}/4\times_{\mathbb{Z}/2}S^1$ is the
union of two circles $C_1$ and $C_2$. Then
$[\Delta]\in\Omega_1(Q_{4n})$ is a nontrivial element.
\end{cor}

\medskip\noindent
{\bf Proof:} All we need, aside from Propositions~\ref{prop:mono}
and \ref{prop:basic}, is the well known fact, \cite{CoFl}
Section~III.20, that the transfer map $t : \Omega_1(H)\rightarrow
\Omega_1(G)$ is defined by $t([M]) = [G\times_H
M]$.\hfill$\square$

\section{Transversality, singular sets and $G$-bordism}
\label{sec:transversality}

The concept of {\em transversality} is central for the whole of
Differential Topology. The basic facts needed in this paper can be
extracted from \cite{GoGu}, Chapter II, \cite{MilnStash}
Chapter~18, or other standard references. As usual, the relation
of transversality is denoted by $f\pitchfork Z$. Here is a list of
basic facts.

\begin{itemize}
\item[${\bf F}_1$]
If $f : M\rightarrow N$ is transverse to a submanifold $Z\subset
N$, then $f^{-1}(Z)$ is a submanifold of $M$ such that ${\rm
codim}_M(f^{-1}(Z)) = {\rm codim}_N(Z)$.
\item[${\bf F}_2$]
Each smooth map $f : M\rightarrow N$ can be perturbed by a
``small'' homotopy to a map $f' : M\rightarrow N$, transverse to
$Z\subset N$. Moreover, if $f$ is already transverse to $Z$ over
an open set $U\subset M$ and $K\subset U$ is a compact set, then
$f'$ can be chosen to coincide with $f$ over $K$.
\item[${\bf F}_3$]
A ``small'' perturbation $f'$ of a map $f : M\rightarrow N$ which
is transverse to $Z$ is also transverse to $Z$.
\item[${\bf F}_4$]
Facts ${\rm F}_1-{\rm F}_3$ together imply that for a given
``arrangement'' ${\cal A}=\{Z_1,\ldots ,Z_k\}$ of submanifolds of
$N$, there always exists a map $f : M\rightarrow N$ such that
$f\pitchfork Z_i$ for all $i=1,\ldots ,k$. In this case we write
$f\pitchfork {\cal A}$.
\item[${\bf F}_5$]
All results remain valid for sections of smooth bundles. In
particular they remain true if $M$ is a free $G$-manifold, where
$G$ is a finite group acting on both $M$ and $N$ as a group of
diffeomorphism and $Z$ is a submanifold of $N$, while all maps are
assumed to be $G$-equivariant. This follows from the fact that
$G$-equivariant maps $f : M\rightarrow N$ are in
$1-1$-correspondence with the sections of the bundle $M\times_G
N\rightarrow M/G$.
\end{itemize}

Suppose that a finite group $G$ acts freely on $S^3$ as a group of
diffeomorphisms. Suppose that $V$ is a real $G$-representation and
${\cal A}$ a $G$-invariant, linear subspace arrangement of $V$. As
usual, let $D({\cal A}):=\cup~{\cal A}$ be the link and $M({\cal
A}):=V\setminus D({\cal A})$ the complement of the arrangement.
Also, let ${\rm Max}({\cal A})= \{P_1,\ldots ,P_s\}$ be the
collection of linear subspaces in ${\cal A}$ of maximum dimension.
Recall that the intersection poset $P = P_{\cal A}$ is an abstract
poset that records the containment relation in ${\cal A}$, i.e.
$(P_{\cal A},\leq)\cong ({\cal A},\subseteq )$. In some papers the
opposite poset $P_{\cal A}^{\rm op}$ is called the intersection
poset of ${\cal A}$. Our choice has the merit that the dimension
function $d : P_{\cal A}\rightarrow N, \, P\mapsto {\rm dim}(P)$
is monotone. Keeping in mind the intended applications in this
paper, we introduce the following, rather restrictive assumption
on the arrangement ${\cal A}$.

\begin{itemize}
\item[${\bf A}_1$]
All maximal elements in ${\cal A}\cong P_{\cal A}$ have dimension
$n-2$ where $n := {\rm dim V}$. Moreover, for any $2$ maximal
elements $P,Q\in P_{\cal A}$, if $P\neq Q$ then ${\rm dim}(P\cap
Q)\leq n-4$.
\end{itemize}

Given a $G$-equivariant map $f : S^3\rightarrow V$, the singular
set $\Delta(f)$ is defined as $\Delta(f) := f^{-1}(D({\cal A}))$.
The singular set $\Delta(f)$ is clearly $G$-invariant. If $f$ is
transverse to all subspaces in the arrangement ${\cal A}$ then, as
a consequence of ${\rm F}_1$, $f^{-1}(C)=\emptyset$ unless ${\rm
dim}(C)=n-2$ i.e. unless $C$ is a maximal subspace in ${\cal A}$.
It follows that $\Delta(f)= \bigcup~\{f^{-1}(C) \mid {\rm
dim}(C)=n-2\}$ is a $1$-dimensional, $G$-submanifold of $S^3$
which has at least $s$ connected components where $s$ is the
cardinality of ${\rm Max}({\cal A})$.

Here are more assumptions on the action of $G$.

\begin{itemize}
\item[${\bf A}_2$]
 The group $G$ preserves the orientation of both
$S^3$ and $V$.
\item[${\bf A}_3$]
Each $P\in {\rm Max}({\cal A})$ is assigned an orientation in such
a way that all these orientations are compatible with the
$G$-action in the sense that the orientation assigned to $g(P)$
agrees with the orientation induced on $g(P)$ from $P\in {\rm
Max}({\cal A})$ by the map $g : P\rightarrow g(P)$.
\end{itemize}

\medskip\noindent
\begin{rem}{\rm
\label{rem:orient} Let ${\rm Stab}(P) = \{g\in G \mid g(P) = P\}$
be the stabilizer of $P\in {\rm Max}({\cal A})$. Than ${\rm
Stab}(P)$ acts on $P$ and it is easy to see that all $Q$ in the
orbit $\{g(P)\}_{g\in G}$ of $P$ can be assigned orientations such
that the condition ${\rm A}_3$ is satisfied if and only if ${\rm
Stab}(P)$ preserves the orientation on $P$. }
\end{rem}

 Under assumptions ${\rm A}_1-{\rm A}_3$, $\Delta(f)$ is an
{\em oriented}, $1$-dimensional, free $G$-manifold which therefore
defines an element $[\Delta(f)]$ in the group $\Omega_1(G)$ of
oriented $G$-bordisms, \cite{CoFl}.

\begin{prop}
 \label{prop:bordism}
Under assumptions ${\rm A}_1-{\rm A}_3$, the element
$[\Delta(f)]\in \Omega_1(G)$ does not depend on the smooth map
$f$. As a consequence if $\Delta = \Delta(f)$ is nontrivial for
{\em one} smooth function $f\pitchfork {\cal A}$, then for any
continuous function $g : S^3\rightarrow V$, the singular set
$\Delta(g)$ is nonempty.
\end{prop}

\medskip\noindent
{\bf Proof:} Suppose that $f$ and $g$ are two maps which are both
transverse to the arrangement ${\cal A}$. Then $\Delta(f)$ and
$\Delta(g)$ are both oriented, $G$-invariant $1$-manifolds and we
are supposed to show that there exists an oriented $G$-bordism
between them. Let $F : S^3\times I\rightarrow V$ be a smooth
homotopy between $f = F(\cdot, 0)$ and $g = F(\cdot, 1)$, for
example $F$ can be obtained by smoothing the linear homotopy
$G(x,t) := (1-t)f(x) + tg(x)$. One can assume that $F$ is
transverse to ${\cal A}$. Indeed, by assumption on $f$ and $g$, it
is already transverse to ${\cal A}$ in a neighborhood $U$ of
$S^3\times\{0,1\}$. By property ${\rm F}_2$, $F$ can be made
transverse to ${\cal A}$ by a small perturbation outside a compact
set $\overline{V\mathstrut}\subset U$ where $V$ is a neighborhood
of $S^3\times\{0,1\}$. Then,
 \[ \Delta(F) = F^{-1}(D({\cal A})) =
 \cup_{P\in {\cal A}}~f^{-1}(P)\subset S^3\times I\]
is an oriented, $G$-invariant, $2$-manifold, which may have
singularities of a very special form. Let us show that these
singularities can be removed and that the desingularized manifold
$\Delta'(F)$ provides the desired $G$-bordism between $\Delta(f)$
and $\Delta(g)$. Let us note that the singular set $S(F)$ has the
following description
 \[ S(F) = \cup\{ F^{-1}(C) \mid {\rm dim}(C) = n-4\}. \]
Note that $S(F)$ is a $G$-invariant subset of $S^3\times I$. Also,
by the property ${\rm F}_1$, the singular set $S(F)$ is
$0$-dimensional. Each $C\in {\cal A}$ of dimension $n-4$ is the
intersection of $2$ or more maximal subspaces in ${\cal A}$ so let
${\cal A}^C := \{P\in {\rm Max}({\cal A}) \mid C\subset P\}$. It
follows that if $x\in F^{-1}(C)$ is a singular point associated to
$C$, then $x$ is a point where manifolds $F^{-1}(P), \, P\in {\cal
A}^C$ intersect transversally. Let $D_{\epsilon}$ be a very small
disc in $F^{-1}(P)$ around $x$. Let us remove the interior
$\stackrel{\circ}{D_{\epsilon}}$ of this disc from $\Delta(F)$ as
well as the interiors of all discs of the form $g(D_{\epsilon})$
for some $g\in G$. Note that the freeness of the action guarantees
that $D_{\epsilon}\cap g(D_{\epsilon})=\emptyset$ for each $g\neq
1$. This means that we enlarged the boundary of $\Delta(F)$ by an
oriented $1$-dimensional, $G$-manifold diffeomorphic to the
manifold $S^1\times G$ with the obvious $G$-action. The process of
removing discs around singular point can be repeated until we
obtain a $2$-dimensional, free $G$-manifold $\Delta_1(F)$ with the
boundary consisting of $\Delta(f), \Delta(g)$ and possibly several
copies of $G\times S^1$. Since $G\times S^1$ is clearly
$G$-cobordant to zero, we can fill the holes with manifolds of the
form $G\times D^2, \, (S^1=\partial(D^2))$ and eventually arrive
at the desired $G$-cobordism $\Delta'(F)$ between $\Delta(f)$ and
$\Delta(g)$. This shows that $\Delta = [\Delta(f)]$ is a well
defined element in $\Omega_1(G)$, independent of $f\pitchfork
{\cal A}$. Suppose that $\Delta\neq 0$ and $\Delta(g) = \emptyset$
for some continuous, $G$-equivariant map $g : S^3\rightarrow V$.
Then there exists a smooth $\epsilon$-approximation $h$ of $g$
such that $h\pitchfork {\cal A}$. If $\epsilon$ is small enough
then $\Delta(h)=\emptyset$ which is in contradiction with $\Delta
= [\Delta(h)]\neq 0$.\hfill$\square$

\section{Partitions by $2$-fans}
\label{sec:2-fans}

\begin{theo}
\label{thm:equi1} Suppose that $n$ and $p$ are two positive
integers such that $n=2k-1$ for some $k$ and $1\leq p\leq n-1$.
Let $\alpha = (\frac{p}{n}, \frac{q}{n})$ where $p+q = n$. Then
there does not exist a $Q_{4n}$-equivariant map $$F : S^3
\rightarrow (W_n)^{\oplus 2}\setminus \cup{\cal A}(\alpha ).$$
\end{theo}

\medskip\noindent
{\bf Proof:} According to Propositions~\ref{prop:spec} and
\ref{prop:two} it is sufficient to prove the theorem in the case
$p=1$. By identifying the space $W:=(W_n)^{\oplus 2}$ with an
appropriate linear subspace of the vector space ${{\rm Mat}}:={\rm
Mat}_{2\times n}(\mathbb{R})$ of all $2$ by $n$ matrices, we may
view ${\cal A}(\alpha)$ as an arrangement in the latter space.
Note that $W$ is $D_{2n}$-invariant, hence $Q_{4n}$-invariant
subspace of ${{\rm Mat}}$, actually $W$ is precisely the
orthogonal complement to the linear space of all $D_{2n}$-fixed
points in ${{\rm Mat}}$. The action of $D_{2n}$ on ${{\rm Mat}}$
is orientation preserving. Indeed, if ${{\rm Mat}}\cong
\mathbb{R}^2\oplus \mathbb{R}^2\oplus\ldots\oplus \mathbb{R}^2 =
\oplus_{i=1}^n~L_i$ is the {\em column decomposition} of ${\rm
Mat}$, then $D_{2n}$ acts essentially by permuting $2$-dimensional
subspaces $L_i$. This shows that the condition ${\rm A}_2$ from
Section~\ref{sec:transversality} is satisfied. Let $\pi_i : {\rm
Mat} \rightarrow \mathbb{R}^2$ be an orthogonal projection on
$L_i\cong \mathbb{R}^2$ and let $S_i := \pi_i^{-1}(0)\subset {\rm
Mat}$. Then, since $\alpha = (\frac{p}{n}, \frac{q}{n}) =
(\frac{1}{n},\frac{n-1}{n})$, the subspaces in ${\cal A}(\alpha)$
of maximum dimension are the spaces $R_i:= W\cap S_i, \,
i=1,\ldots, n$. It follows that the condition ${\rm A}_1$ is also
satisfied. Finally, since ${\rm Stab}(R_1) = \{I,JE\}\cong
\mathbb{Z}/2$ and $JE$ fixes all vectors in $R_1$, the condition
${\rm A}_3$ in Section~\ref{sec:transversality} and
Proposition~\ref{prop:bordism} is also fulfilled.

In order to apply Proposition~\ref{prop:bordism}, one is supposed
to select carefully a $Q_{4n}$-equivariant map $f: S^3\rightarrow
W$ transverse to the arrangement ${\cal A}(\alpha)$, compute the
associated singular set $\Delta(f)$ and show that the
corresponding element $[\Delta(f)]\in\Omega_1(Q_{4n})$ is
nontrivial.

\medskip
In the construction of the map $f$, we will use both the smooth
and a simplicial model for the sphere $S^3$. As a complex vector
space, the quaternions have the decomposition $\mathbb{H}\cong
\mathbb{C}_{(1)}\oplus \mathbb{C}_{(2)}$ where $1,i\in
\mathbb{C}_{(1)}$ and $j,k\in \mathbb{C}_{(2)}$. Each unit
quaternion $q\in S^3$ has a unique decomposition of the form $q =
(\cos\alpha) q_1 + (\sin\alpha) q_2$ where $0\leq
\alpha\leq\frac{\pi}{2}$ and $q_1$ and $q_2$ are unit quaternions
in $\mathbb{C}_{(1)}$ and $\mathbb{C}_{(2)}$ respectively. This
simple observation is at the root of the well known {\em join}
decomposition $S^3\cong S^1_{(1)}\ast S^1_{(2)}$ where $S^1_{(1)}$
and $S^1_{(2)}$ are respectively the unit circles in
$\mathbb{C}_{(1)}$ and $\mathbb{C}_{(2)}$. Let $P_{(1)}$ and
$P_{(2)}$ be regular polygons in $\mathbb{C}_{(1)}$ and
$\mathbb{C}_{(2)}$ with respective vertices $a_p :=
\cos\frac{p\pi}{n} + i\sin\frac{p\pi}{n}$ and $b_q :=
j\cos\frac{q\pi}{n} - k \sin\frac{q\pi}{n} k$ for $p,q=0,\ldots
,2n-1$. This enumeration of vertices of the polygon $P_{(2)}= {\rm
conv}\{b_q\}_{q=1}^n$ is chosen because we work with left actions
and we want to have $b_q = ja_q$. The (simplicial) join $P:=
P_{(1)}\ast P_{(2)}$ is a triangulation of the sphere $S^3$. We
shall construct a simplicial  map $f : P \rightarrow W$ which is
$Q_{4n}$-equivariant. Although strictly speaking $f$ is not a
smooth map, $\Delta(f)$ still provides the information needed for
an application of Proposition~\ref{prop:bordism}. Indeed, one can
from the start interpret $P$ as a smooth triangulation of $S^3$ or
alternatively recall that there is no essential difference between
the smooth and $PL$ category in dimension $3$ so in particular all
results in Section~\ref{sec:transversality} remain valid for
$PL$-manifolds.

\medskip
Each $L_i\cong \mathbb{R}^2_{(i)}$ in the column decomposition
${\rm Mat} = \oplus_{i=1}^n~L_i$ is identified with the complex
plane $\mathbb{C} = \mathbb{C}_{(i)}$. Choose a nonzero vector
$v_0\in L_1 = \mathbb{C}_{(1)}$. Then $v_i := \epsilon^{2i}v_0$
are vertices of a regular $n$-gon $E_n$ in $L_1 =
\mathbb{C}_{(1)}$.

\medskip\noindent
{\bf Caveat:} As before we allow the indices $p, q$ etc. of
vertices $a_p, b_q$ of the complex $P$ to range over all integers
and assume that $x_p = x_{p'}$ whenever  $p =_{2n} p'$ are
congruent modulo $2n$. Similarly, the indices $j$ of vertices
$v_j$ also range over $\mathbb{Z}$ with the convention that $v_i =
v_j$ if $i =_n j$.

\medskip
We want to have $f(\epsilon a_0)= Ef(a_0) = \epsilon^2\cdot
f(a_0)$ so $f(a_0)$ is defined by
 \[
 f(a_0) := v_0\oplus\epsilon^2 v_0\oplus\ldots\oplus
\epsilon^{2n-2} v_0 = v_0\oplus v_1\oplus\ldots\oplus v_{n-1}\in
\oplus_{i=1}^n~L_i .
 \]
 The condition $f(g\cdot a_0) = g\cdot
f(a_0)$ forces us to define $f(a_p)=f(\epsilon^p a_0)$ and $f(b_q)
= f(j\cdot a_q)$ by
 \begin{equation}
 \label{eqn:def1}
f(a_p) = \epsilon^{2p} v_0\oplus\epsilon^{2p+2}
v_0\oplus\ldots\oplus \epsilon^{2p+2n-2} v_0 = v_p\oplus
v_{p+1}\oplus\ldots\oplus v_{p+n-1}
 \end{equation}
\begin{equation}
 \label{eqn:def2}
f(b_q)= JE^{q}(v_0\oplus v_1\oplus\ldots\oplus v_{n-1}) =
v_{q+n-1}\oplus v_{q+n-2}\oplus\ldots\oplus v_{q}
 \end{equation}
so for example,
\[
f(b_0) = f(j\cdot a_0) = v_{n-1}\oplus\ v_{n-2}\oplus\ldots\oplus
v_0.
 \]

The map $f$, already  defined as a $Q_{4n}$-equivariant map on the
set of vertices of $P = P_{(1)}\ast P_{(2)}$, admits a unique
simplicial extension on the whole of $P$. Let us show that $f$ is
$Q_{4n}$-equivariant.  This will follow once we convince ourselves
that $P = P_{(1)}\ast P_{(2)}$ is a free $Q_{4n}$-complex. Note
that a $3$-simplex $\sigma = a_pa_{p+1}b_qb_{q+1}\in P$ is
uniquely reconstructed from its label $l(\sigma) = (p,q)\in
\mathbb{Z}/2n\times \mathbb{Z}/2n$. The induced action of $Q_{4n}$
on the set $\mathbb{Z}/2n\times \mathbb{Z}/2n$ of labels of all
$3$-simplices is obviously free since $\epsilon(p,q) = (p+1,q-1)$
and $j(p,q) = (q+n,p)$. Indeed if for example
$j\epsilon^m(p,q)=(p,q)$, then $(q-m+n,p+m)=(p,q)$ which implies
$p=_{2n}p+n$, a contradiction.

\medskip
Let us determine the singular set $\Delta(f)\subset P$ of $f$,
proving along the way that $f$ is transverse to the arrangement
${\cal A}$. The maximal subspaces in the arrangement ${\cal A}$
are $R_i = {\rm Ker}(\pi_i)\cap W$ where $\pi_i : {\rm
Mat}\rightarrow L_i$ is the projection. Then,

\begin{equation}
\label{eqn:delta}
 \Delta(f) = \bigcup_{i=1}^n f^{-1}(R_i) = \bigcup_{i=1}^n (\pi_i\circ
 f)^{-1}(0) .
\end{equation}
It is sufficient to determine $O := (\pi_1\circ f)^{-1}(0)$ as a
(left) $\mathbb{Z}/4$-manifold. Indeed, the stabilizer of $R_1$ in
$Q_{4n}$ is the group
 \[
 {\rm Stab}(R_1) = \{1,j\epsilon, -1,-j\epsilon\}\cong
 \mathbb{Z}/4.
 \]
Then $O$ is a ${\rm Stab}(R_1)$-manifold and because of
(\ref{eqn:delta}), $\Delta(f)\cong Q_{4n}\times_{\mathbb{Z}/4} O$
as a (left) $Q_{4n}$-manifold. Note that $\pi_i f(a_p) =
v_{p+i-1}$ and $\pi_i f(b_q) = v_{q+n-i}$ for all $i=1,\ldots ,n$.
Hence, for $i=1$
 \[
 \pi_1 f(a_p) = v_p      \qquad \mbox{ {\rm and } } \qquad
 \pi_1 f(b_q) = v_{q+n-1}.
 \]

A triangle $a_pa_{p+1}b_q\in P = P_{(1)}\ast P_{(2)}$ is called
{\em good} if it intersects $O = (\pi_1\circ f)^{-1}(0)$ i.e. if
and only if
\begin{equation}
\label{eqn:cond1}
 0\in {\rm conv}\{\pi_1f(a_p),\pi_1f(a_{p+1}),\pi_1f(b_q)\} =
 {\rm conv}\{v_p, v_{p+1}, v_{q+n-1}\}.
\end{equation}
Similarly, a triangle $a_pb_qb_{q+1} \in P$ is {\em good} if it
intersects $O = (\pi_1\circ f)^{-1}(0)$ which happens if and only
if
\begin{equation}
\label{eqn:cond2}
 0\in {\rm conv}\{\pi_1f(a_p),\pi_1f(b_{q}),\pi_1f(b_{q+1})\} =
 {\rm conv}\{v_p, v_{q+n-1}, v_{q+n}\}.
\end{equation}
Note that if $a_pa_{p+1}b_q$ is a good triangle, then the same
property have the triangles
\begin{equation}
\label{eqn:temena} \ldots \quad a_{p}b_{q-1}b_q \quad
a_pa_{p+1}b_q \quad a_{p+1}b_qb_{q+1} \quad a_{p+1}a_{p+2}b_{q+1}
\quad a_{p+2}b_{q+1}b_{q+2} \quad \ldots
\end{equation}
Let us show that the $1$-manifold $O$ is actually the union of two
circles, $O = C_1\cup C_1$. The triangle $a_0a_1b_{k}$, where $n =
2k-1$, intersects one of these circles, say $C_1$ while the
triangle $a_0a_1b_{3k-1} = a_0a_1b_{n+k}$ intersects the other.
Actually the circles, connected components of $O$, are in $1$--$1$
correspondence with chains of consecutive triangles of the form
(\ref{eqn:temena}). It is than easy to check that there are
precisely two such chains. Alternatively, one can show that both
$a_pa_{p+1}b_q$ and $a_{p+1}b_{q}b_{q+1}$ are good triangles if
and only if either $q-p = k$ or $q-p = n+k = 3k-1$. Then, if both
$a_pa_{p+1}b_q$ and $a_{p+1}b_{q}b_{q+1}$ are labelled by
$(p,q)\in \mathbb{Z}/2n\times \mathbb{Z}/2n$, then the set
$L_{\Delta}$ of labels of all good triangles is described as
 \begin{equation}
 \label{eqn:labels}
 L_{\Delta}= \theta^{-1}(k)\cup \theta^{-1}(3k-1)
 \end{equation}
where $\theta : \mathbb{Z}/2n\times \mathbb{Z}/2n \rightarrow
\mathbb{Z}/2n$ is the map defined by $\theta(p,q) = q-p$. In this
representation the sets of labels $\theta^{-1}(k)$ and
$\theta^{-1}(3k-1)$ correspond to circles $C_1$ and $C_2$
respectively. The meaning of labels becomes even more transparent
if we observe that the triangles $a_pa_{p+1}b_q$ and
$a_{p+1}b_{q}b_{q+1}$ are faces of a common $3$-simplex
$\sigma_{(p,q)}= a_pa_{p+1}b_qb_{q+1}\in P$. Note also that if
$\gamma : \mathbb{Z}/2n\times \mathbb{Z}/2n \rightarrow
\mathbb{Z}/n$ is the map defined by $\gamma(p,q) =
[\theta(p,q)]_{\mbox{$\scriptstyle\rm mod\,  n$}}$, then the
decomposition $\Delta(f) = \cup_{i=1}^n~\gamma^{-1}(i)$
corresponds to the decomposition (\ref{eqn:delta}) of the singular
set $\Delta(f)$.

A consequence of this analysis is that $\Delta(f)$ is, as an
oriented $Q_{4n}$-manifold, isomorphic to
$Q_{4n}\times_{\mathbb{Z}/2}S^1$. Hence, by
Corollary~\ref{cor:vazan}, $[\Delta(f)]$ is a nontrivial element
$\Delta$ of $\Omega_1(Q_{4n})$. By Proposition~\ref{prop:bordism},
$\Delta(g)$ is nonempty for each continuous, $Q_{4n}$-equivariant
map $g : S^3\rightarrow (W_n)^{\oplus 2}$, i.e. there does not
exist a $Q_{4n}$-equivariant map $g : S^3\rightarrow (W_n)^{\oplus
2}\setminus {\cal A}(\alpha)$. This completes the proof of the
theorem.\hfill$\square$

\begin{theo}
\label{thm:prva} Suppose that $\alpha = (s,t)$ is vector in
$\mathbb{R}^2$ such that $s+t=1$ and $s,t>0$. Then any collection
of three proper measures $\mu_1,\mu_2,\mu_3$ on the sphere $S^2$
admits an $\alpha$-partition by a $2$-fan ${\mathfrak p} = (x;
l_1,l_2)$. In other words, \[{\cal A}_{3,2}=\{(s,t)\in
\mathbb{R}^2 \mid s+t =1 { \mbox{ \rm and } } s,t > 0\}.\]
\end{theo}

\medskip\noindent
{\bf Proof:} A simple limit and compactness argument shows that
${\cal A}_{3,2}$ is a closed subspace of $T := \{(s,t)\in
\mathbb{R}^2 \mid s+t =1 \mbox{ \rm and } s,t>0\}$. Hence it is
sufficient to show that $S\subset{\cal A}_{3,2}$ for a dense
subset $S\subset T$. Let $(s,t) = (\frac{p}{n}, \frac{q}{n})$,
where $n$ is an odd number. In order to show that $(\frac{p}{n},
\frac{q}{n})\in {\cal A}_{3,2}$ it is sufficient, according to
Proposition~\ref{prop:reduction}, to show that there does not
exist a $D_{2n}$-equivariant map $f : V_2(\mathbb{R}^3)\rightarrow
(W_n)^{\oplus 2}\setminus {\cal A}(\alpha)$. In light of
Proposition~\ref{prop:homological} applied on $X=S^3, G=Q_{4n}$
and $H= \mathbb{Z}/2$, this is precisely the statement of
Theorem~\ref{thm:equi1}.\hfill$\square$

\section{Partitions by $3$-fans}
\label{sec:3-fans}

The central part of the proof of Theorem~\ref{thm:equi1} was the
construction of a special, $Q_{4n}$-equivariant map $f :
S^3\rightarrow W_n$ and the evaluation of its singular set
$\Delta(f)$. A similar strategy is applied in the case of
partitions by $3$-fans. Unfortunately, the invariant $\Delta =
[\Delta(f)]\in\Omega_1(Q_{4n})$ turns out to be zero in this case,
Proposition~\ref{prop:vanishes}. We nevertheless reproduce here
this computation for two reasons. The combinatorial details appear
to be interesting in themselves and, since the first obstruction
to the existence of an equivariant map vanishes, we expect that
they may be useful in the computation of the secondary
obstruction. Secondly, the evaluation of $\Delta$ in the case of
$3$-partitions provides a neat example what in principle can ``go
wrong'' with this approach.

We keep the same notation as before so in particular $\mathbb{H}
\cong \mathbb{C}_{(1)}\oplus \mathbb{C}_{(2)}$ and $S^3 =
S^1_{(1)}\ast S^1_{(2)}$ while $P = P_{(1)}\ast P_{(2)}$ is the
triangulation of the $3$-sphere obtained as the join of regular
$(2n)$-gons $P_{(1)}$ and $P_{(2)}$ in $\mathbb{C}_{(1)}$ and
$\mathbb{C}_{(2)}$ respectively. The indexing of vertices
$\{a_p\}_{p=1}^{2n}$ of $P_{(1)}$ is also unchanged so as before
$a_p = \epsilon^p a_0 = \cos\frac{p\pi}{n} + i\sin\frac{p\pi}{n}$.
However, the enumeration of vertices $\{b_q\}_{q=1}^{2n}$ of
$P_{(2)} = j(P_{(1)})$ is different. In order to make the
calculations more transparent we assume that $j a_p = b_q$ is
equivalent to $p+q = n+1$ (mod $2n$) or in other words $j a_p =
b_{n-p+1}$ for each $p = 0,1,\ldots, 2n-1$. From here one easily
deduces the following equalities
\begin{equation}
\label{eqn:tablica}
\begin{array}{lccl}
 j a_p = b_{n-p+1} & & & j b_q = a_{-q+1} \\
 \epsilon a_p = a_{p+1} & & &  \epsilon b_q = b_{q+1}.
\end{array}
\end{equation}
Let $v_i = e_i - (1/n)\sum_{j+1}^n e_j, \, i=1,\ldots, n$ be the
vertices of the regular, $\mathbb{D}_{2n}$-invariant simplex in
$W_n$. Let $f' : P^{(0)}\rightarrow W_n$ be the map defined on the
$0$-skeleton $P^{(0)} = \{a_p\}_{p=1}^{2n}\cup\{b_q\}_{q=1}^{2n}$
of $P$ by
\begin{equation}
\label{eqn:map} f'(a_p) = v_p \qquad \mbox{  \rm and  } \qquad
f'(b_q) = v_q .
\end{equation}
Let $f : P\rightarrow W_n$ be the simplicial extension of $f'$ on
the complex $P = P_{(1)}\ast P_{(2)}\cong S^3$. Then,
\begin{equation}
\label{eqn:equitest}
\begin{array}{ccccccccc}
f(ja_p) &=& f(b_{n-p+1}) &=& v_{n-p+1} &=& Jv_p &=& Jf(a_p)\\
f(\epsilon a_p) &=& f(a_{p+1}) &=& v_{p+1} &=& E v_p &=& Ef(a_p)
\end{array}
\end{equation}
which means that $f'$ is $\mathbb{Q}_{4n}$-equivariant. As in the
proof of Theorem~\ref{thm:prva}, $f'$ is extended to a unique
$\mathbb{Q}_{4n}$-equivariant, simplicial map \[ f : P_{(1)}\ast
P_{(2)}\rightarrow W_n. \]
 The next step is to determine the associated singular set
\begin{equation}
\label{eqn:determine}
\begin{array}{ccccccc}
\Delta(f) &=& f^{-1}(\cup{\cal A}(\alpha)) &=& \bigcup_{g\in
\mathbb{D}_{2n}} f^{-1}(gL(\alpha)) &=& \bigcup_{h\in
\mathbb{Q}_{4n}} hf^{-1}(L(\alpha)).
\end{array}
\end{equation}
Let us write the defining equations (\ref{eqn:L}) and
(\ref{eqn:formsforms}) for $L(\alpha)$ in the special case $\alpha
= (\frac{p}{n}, \frac{q}{n}, \frac{r}{n})$. If $P:= p$ and $Q:=
p+q$ then the equations $z_1 = z_2 = z_3 = 0$ can be rewritten as
\begin{equation} \label{eqn:reduce}  x_1+\ldots + x_P  =
x_{P+1}\ldots + x_Q = x_{Q+1}+\ldots +x_n = 0.
\end{equation}
In light of the condition $z_1+z_2+z_3 = 0$, one of the equations
in (\ref{eqn:reduce}) is redundant so $L(\alpha)$ is a subspace of
$W_n$ of codimension $2$. Clearly $\Delta(f) \cong
Q_{4n}\times_H~O$, where $O := f^{-1}(L(\alpha))$ and $H:= {\rm
Stab}(L(\alpha))\subset Q_{4n}$ is the stabilizer of $L(\alpha)$.
It is clear that the group $H = H_{\alpha}$ depends on the vector
$\alpha = (\frac{p}{n}, \frac{q}{n}, \frac{r}{n})$. More
precisely, $H_\alpha = \{-1,1\}$ if all three integers $p,q,r$ are
pairwise distinct, $H_\alpha = \{1, j\epsilon^q, -1,
-j\epsilon^q\}\cong\mathbb{Z}/4$ if $\alpha = (\frac{p}{n},
\frac{p}{n}, \frac{q}{n})$ and $p\neq q$, and $H_\alpha \cong S_3$
if $p=q=r$.

It turns out that the $H$-manifold $O$ is always the union of $2$
circles, $O = O_1\cup O_2$. This assertion is proved by a careful
analysis of which simplices $a_pa_{p+1}b_qb_{q+1}$ intersect
$\Delta(f)$ and $O$. The triangles $a_pa_{p+1}b_q$ and
$a_pb_qb_{q+1}$ are called {\em good} if they intersect $O =
f^{-1}(L(\alpha))$. Let $[1,2n] = {\cal P}\cup{\cal Q}\cup{\cal
R}$ be the partition of the index set $[1,2n] = \{1,2,\ldots
,2n\}$ defined by ${\cal P} = [1,P]\cup [n,P+n], {\cal Q} =
[P+1,Q]\cup [P+n+1,Q+n], {\cal R} = [Q+1,n]\cup [Q+n+1,2n].$

\medskip\noindent
{\bf Claim~1:} A triangle $a_pa_{p+1}b_q$ is {\em good} if and
only if $\{p,p+1,q\}$ intersects each of the sets ${\cal P}, {\cal
Q}$ and ${\cal R}$. Similarly, a triangle $a_pb_qb_{q+1}$ is {\em
good} if and only if $\{p,q,q+1\}$ has a nonempty intersection
with each of the sets ${\cal P}, {\cal Q}$ and ${\cal R}$.

\placedrawing[htbp]{lim2.lp}{}{lim2}

\medskip\noindent
An example of a good triangle is $a_{2n}a_1b_{P+1}$ which
intersects $O$ precisely in the point $c =
(r/n)a_{2n}+(p/n)a_1+(q/n)b_{P+1}\in P$. More generally, it is
easily verified that if $x_iy_jz_k$ is a good triangle such that
$i\in {\cal P}, j\in {\cal Q}$ and $k\in {\cal R}$, then $c =
(p/n)x_i+(q/n)y_j+(r/n)z_k$ is the unique point in  $x_iy_jz_k$
which intersects $O$. This observation is essentially all we need
for the proof of Claim~1 and at the same time it shows that the
map $f$ is transverse to $L(\alpha)$. Two good triangles are
called {\em adjacent} if they are faces of a common $3$-simplex in
$P$. They are in the same component  if they can be connected by a
chain of adjacent triangles. It turns out that the set $\Gamma$ of
all good triangles falls apart into $2$ components, $\Gamma =
\Gamma_1\cup \Gamma_2$. The first component $\Gamma_1$, the one
that contains the triangle $a_{2n}a_1b_{P+1}$, is partially
tabulated by (\ref{eqn:listgood}). The list for the second
component $\Gamma_2$ is obtained from by (\ref{eqn:listgood}) by
formally interchanging the letters $a$ and $b$. For example
$a_{2n}a_1b_{Q+n}$ and $a_{P+1}b_{2n}b_1$ are in $\Gamma_2$, while
$a_{Q+n}b_{2n}b_1$ and $a_{2n}a_1b_{P+1}$ are in $\Gamma_1$.

\begin{equation}
\label{eqn:listgood}
\begin{array}{lllllll}
{\scriptstyle a_{2n}a_1b_{P+1}} && {\scriptstyle a_{2n}a_1b_{P+2}}
&& \ldots && {\scriptstyle  a_{2n}a_1b_Q}\\ {\scriptstyle
a_{1}b_Qb_{Q+1}} && {\scriptstyle a_2b_Qb_{Q+1}} && \ldots &&
{\scriptstyle a_Pb_Qb_{Q+1}}\\ {\scriptstyle a_{P}a_{P+1}b_{Q+1}}
&& {\scriptstyle a_{P}a_{P+1}b_{Q+2}} && \ldots && {\scriptstyle
a_{P}a_{P+1}b_n}\\ {\scriptstyle a_{P+1}b_nb_{n+1}} &&
{\scriptstyle a_{P+1}b_nb_{n+1}} && \ldots && {\scriptstyle
a_Qb_nb_{n+1}}\\ {\scriptstyle a_{Q}a_{Q+1}b_{n+1}} &&
{\scriptstyle a_{Q}a_{Q+1}b_{n+2}} && \ldots && {\scriptstyle
a_{Q}a_{Q+1}b_{P+n}}\\ {\scriptstyle a_{Q+1}b_{P+n}b_{P+n+1}} &&
{\scriptstyle a_{Q+2}b_{P+n}b_{P+n+1}} && \ldots && {\scriptstyle
a_nb_{P+n}b_{P+n+1}}\\ {\scriptstyle a_{n}a_{n+1}b_{P+n+1}} &&
\ldots && \ldots && \ldots\\ \ldots &&\ldots &&\ldots &&\ldots
\end{array}
\end{equation}

\begin{prop}
\label{prop:vanishes} Suppose that $n$ is odd, $\alpha =
(\frac{p}{n}, \frac{q}{n}, \frac{r}{n})$ and the integers $p,q,r$
are pairwise distinct. Then the obstruction element $\Delta =
[\Delta(f)]\in\Omega_1(Q_{4n})$ is zero, where $\Delta(f)$ is the
singular set of the map $f : S^3\rightarrow W_n$ described by {\rm
(\ref{eqn:map})} and {\rm (\ref{eqn:equitest})}. Moreover, under
these conditions the first obstruction for the existence of a
$Q_{4n}$-equivariant map $f : S^3\rightarrow W_n\setminus\cup{\cal
A}(\alpha)$ also vanishes.
\end{prop}

\medskip\noindent
{\bf Proof:} The assumptions on $n$ and $\alpha$ guarantee that
the conditions ${\rm A}_1-{\rm A}_3$ of
Proposition~\ref{prop:bordism} are satisfied. The analysis of the
singular set $\Delta(f)$ reveals that $\Delta(f)\cong Q_{4n}\times
S^1$ as a $Q_{4n}$-set, hence $\Delta =
[\Delta(f)]\in\Omega_1(Q_{4n})$ is a trivial element.

The first obstruction to the existence of an equivariant map $g :
S^3\rightarrow W_n\setminus{\cal A}(\alpha)$ is an element
$\omega\in H^2_{Q_{4n}}(S^3,\pi_1(W_n\setminus{\cal A}(\alpha)))$.
Condition ${\rm A}_1$ implies that
\[
A:= \pi_1(W_n\setminus{\cal A}(\alpha))\cong H_1(W_n\setminus{\cal
A}(\alpha); \mathbb{Z})\cong\mathbb{Z}^{2n}.
\]
As a $D_{2n}$-module, $A \cong D_{2n}\times\mathbb{Z}$ is
isomorphic to the regular representation of $D_{2n}$ over
$\mathbb{Z}$. This means that as a $Q_{4n}$-module, $A\cong
Q_{4n}\times_{\mathbb{Z}/2}\mathbb{Z}$ for an appropriate subgroup
$\mathbb{Z}/2\subset Q_{4n}$. It follows that $H^2_{Q_{4n}}(S^3,
A)\cong H^2_{\mathbb{Z}/2}(S^3,\mathbb{Z})$ by the ``extensions of
scalars'' isomorphism, \cite{Brown}. By eqivariant Poincar\' e
duality, \cite{Wall}, $H^2_{\mathbb{Z}/2}(S^3,\mathbb{Z}^2)\cong
H_1^{\mathbb{Z}/2}(S^3, \mathbb{Z})$ and the obstruction class
$\omega$ corresponds to the class $\hat{\omega}\in
H_1^{\mathbb{Z}/2}(S^3, \mathbb{Z}^2)$ represented by two circles.
Finally, the triviality of $\hat{\omega}$ follows from the
geometric fact that two circles represent a trivial element in
$$
H_1^{\mathbb{Z}/2}(S^3, \mathbb{Z})\cong  H_1^{\mathbb{Z}/2}(S^3,
\mathbb{Z})\cong \Omega_1(\mathbb{Z}/2). \eqno{\square}
$$


\vfill\newpage

\begin{thebibliography}{}

\bibitem{Aki2000}
J.~Akiyama, A.~Kaneko, M.~Kano, G.~Nakamura, E.~Rivera-Campo,
S.~Tokunaga, and J.~Urrutia.
\newblock Radial perfect partitions of convex sets in the plane.
\newblock In {\em Discrete and Computational Geometry} (J.~Akiyama
et al. eds.), {\em Lect. Notes Comput. Sci. 1763, pp. 1--13.}
Springer, Berlin 2000.

\bibitem{Bar93}
I.~B\' ar\' any,
\newblock Geometric and combinatorial applications of Borsuk's theorem,
\newblock {\it New trends in Discrete and Computational Geometry},
J\' anos Pach, ed.,
\newblock Algorithms and Combinatorics 10,
Springer-Verlag, Berlin, 1993.

\bibitem{BaMa2001}
I.~B\' ar\' any, J.~Matou\v sek.
\newblock Simultaneous partitions of measures by $k$-fans,
\newblock {\em Discrete Comp. Geometry}, 25\, (2001), 317--334.

\bibitem{BaMa2002}
I.~B\' ar\' any, J.~Matou\v sek.
\newblock Equipartitions of two measures by a $4$-fan,
(preprint).

\bibitem{Bes2000}
S.~Bespamyatnikh, D.~Kirkpatrick, and J.~Snoeyink.
\newblock Generalizing ham sandwich cuts to equitable
subdivisions.
\newblock {\em Discrete Comput. Geom.,} 24:605--622, 2000.

\bibitem{Bjo91}
A.~Bj{\"o}rner.
\newblock Topological methods,
\newblock In R.~Graham, M.~Gr{\"o}tschel, and L.~Lov{\'a}sz, editors,
{\it   Handbook of Combinatorics}. North-Holland, Amsterdam, 1995.

\bibitem{Brown}
K.S.~Brown.
\newblock {\it Cohomology of groups},
\newblock Springer-Verlag, New York, Berlin, 1982.

\bibitem{CaEi}
H.~Cartan and S.~Eilenberg.
\newblock {\it Homological Algebra}
\newblock Princeton University Press, 1956.

\bibitem{CoFl}
P.E.~Conner and E.E.~Floyd.
\newblock {\it Differentiable periodic maps},
Springer-Verla, Berlin 1964.

\bibitem{Dieck87}
T.~tom~Dieck,
\newblock {\it Transformation groups},
\newblock de Gruyter Studies in Math. 8,
Berlin, 1987.


\bibitem{GoGu}
M.~Golubitsky, V.~Guillemin.
\newblock {\em Stable Mappings and Their Singularities}.
\newblock Graduate Texts in Mathematics 14, Springer--Verlag 1973.

\bibitem{Go}
C.M.~Gordon.
\newblock The $G$-signature theorem in dimension $4$,
\newblock in {\it A la recherche de la topolgie perdue},
L.~Guillou, A.~Marin (eds.), Birkh\" auser, 1986.

\bibitem{Ito2000}
H.~Ito, H.~Uehara, and M.~Yokoyama.
\newblock $2$-dimension ham-sandwich theorem for partitioning into
three convex pieces.
\newblock In {\em Discrete and Computational Geometry} (J.~Akiyama
et. al eds.), {\em Lect. Notes Comput. Sci. 1763}, pp. 129--157.
Springer, Berlin 2000.

\bibitem{Kan1999}
A.~Kaneko, M.~Kano.
\newblock Balanced partitions of two sets of points in the plane.
\newblock {\em Comput. Geom. Theor. Appl.}, 13(4), 253--261, 1999.

\bibitem{MaVrZi}
P.~Mani-Levitzka, S.~Vre\' cica, R.~\v Zivaljevi\' c,
\newblock Combinatorics and topology of partitions of masses
by hyperplanes,
\newblock (in preparation).

\bibitem{Mat94}
J.~Matou\v sek,
\newblock Topological methods in Combinatorics and Geometry,
\newblock {\it Lecture notes, Prague  1994.}
(updated version, February 2002,
www.ms.mff.cuni.cz/~matousek/lecturenotes.html).


\bibitem{MilnStash}
J.W.~Milnor, J.D.~Stasheff.
\newblock {\em Characteristic Classes}.
\newblock Annals of Mathematics Studies 76, Princeton University
Press 1974.


\bibitem{OrlikTerao}
P.~Orlik, H.~Terao.
\newblock {\em Arrangements of Hyperplanes.}
\newblock Grundlehren der mathematischen Wissenschaften 300,
Springer-Verlag 1992.



\bibitem{Pach93}
J.~Pach (Ed.),
\newblock {\it New Trends in Discrete and Computational Geometry},
\newblock Algorithms and Combinatorics 10,
 Springer 1993.


\bibitem{Ramos}
E.~Ramos,
\newblock Equipartitions of mass distributions,
by hyperplanes.
\newblock {\it Discrete  Comput. Geom.,} 15 : 147--167, 1996.

\bibitem{Sak1998}
T.~Sakai.
\newblock Radial partitions pf point sets in $R^2$.
\newblock Manuscript, Tokoha Gakuen University, 1998.


\bibitem{TV93}
H. Tverberg, S. Vre\' cica,
\newblock On generalizations of Radon's theorem and the
ham sandwich theorem, \newblock {\it Europ. J. Combinatorics} 14,
1993, pp. 259--264.


\bibitem{VZ92}
S. Vre\' cica, R. \v Zivaljevi\' c,
\newblock The ham sandwich theorem revisited,
\newblock {\it Israel J. Math.}
78,  1992, pp. 21--32.

\bibitem{VreZiv2001}
S. Vre\' cica, R. \v Zivaljevi\' c.
\newblock Conical equipartitions of massdistributions.
\newblock {\em Discrete Comput. Geom.}, 225:335--350, 2001.

\bibitem{VreZiv2002}
S. Vre\' cica, R. \v Zivaljevi\' c.
\newblock  Arrangements, equivariant maps and partitions of measures
by $4$-fans,
\newblock (preprint).

\bibitem{Wall}
C.T.C. Wall.
\newblock {\em Surgery on Compact Manifolds}.
\newblock Academic Press, 1970.


\bibitem{elisat}
R.T.~\v Zivaljevi\' c,
\newblock Topological methods,
\newblock in  {\it CRC Handbook of
Discrete and Computational Geometry}, J.E.~Goodman, J.~O'Rourke,
eds.
\newblock CRC Press, Boca Raton 1997.


\bibitem{T-V}
R.T.~\v Zivaljevi\' c,
\newblock The Tverberg--Vre\' cica problem and the combinatorial
geometry on vector bundles,
\newblock {\it Israel J. Math.} 111 (1999), 53--76.

\bibitem{guide1}
R.T.~\v Zivaljevi\' c,
\newblock User's guide to equivariant methods in combinatorics,
\newblock {\it Publ. Inst. Math. Belgrade}, 59(73), 1996, 114--130.

\bibitem{guide2}
R.T.~\v Zivaljevi\' c,
\newblock User's guide to equivariant methods in combinatorics II,
\newblock {\it Publ. Inst. Math. Belgrade}, 64(78), 1998, 107--132.

\bibitem{ZV90}
R.~\v Zivaljevi\' c, S.~Vre\' cica,
\newblock An extension of the ham sandwich theorem,
\newblock {\it Bull. London Math. Soc.} 22,  1990, pp. 183--186.



\end{thebibliography}
\end{document}